\documentclass[a4paper,11pt]{article}
\usepackage{url}
\usepackage{latexsym}
\usepackage{amssymb}
\usepackage{amsmath}
\usepackage{amsthm}
\usepackage{times}
\usepackage{xspace}
\usepackage[dvips]{graphicx}

{\bfseries}{\itshape}
{\bfseries}{\itshape}
{\bfseries}{\itshape}
{\bfseries}{\itshape}
\newtheorem{theorem}{Theorem}{\bfseries}{\itshape}
\newtheorem{lemma}{Lemma}{\bfseries}{\itshape}
\newtheorem{definition}{Definition}{\bfseries}{\itshape}

\title{Boxicity of Series Parallel Graphs}
\author{Ankur Bohra
\thanks{Dept. of Mathematics, Indian Institute of Technology, Delhi, 110016.
Email : {\tt mau02481@ccsun50.iitd.ernet.in}}
\and L. Sunil Chandran
\thanks{Dept. of Computer Science and Automation, Indian Institute of Science, Bangalore, 560012. (\textbf{Corresponding Author})Email : 
{\tt sunil@csa.iisc.ernet.in}}
\and J. Krishnam Raju 
\thanks{Dept. of Computer Science and Automation, Indian Institute of Science, Bangalore, 560012. Email : {\tt jkraju@csa.iisc.ernet.in}}
}
\date{}

\begin{document}
\pagestyle{plain}
\pagenumbering{arabic}
\maketitle

\newcommand{\mr}{\ensuremath{\mathbf{\eta}}}
\newcommand{\len}{\ensuremath{t}}
\newcommand{\minor}{\ensuremath{\preceq}}
\newcommand{\node}[1]{\ensuremath{\langle #1 \rangle }}
\newcommand{\floor}[1]{\ensuremath{\lfloor #1 \rfloor }}
\newcommand{\ceil}[1]{\ensuremath{\lceil #1 \rceil }}
\newcommand{\cart}{\ensuremath{\,\Box\,}}
\newcommand{\brac}[1]{\ensuremath{\lbrace #1 \rbrace }}
\newcommand{\sbrac}[1]{\ensuremath{\lbrack #1 \rbrack }}
\newcommand{\ignore}[1]{}

\begin{abstract}
The three well-known graph classes, planar graphs(${\cal P}$), series-parallel graphs (${\cal SP}$) and outer
planar graphs(${\cal OP}$) satisfy the following proper inclusion relation: ${\cal OP} \subset {\cal SP} \subset {\cal P}$. 
It is known that $\mbox{box}(G) \le 3$ if $G \in {\cal P}$ and $\mbox{box}(G) \le 2$ 
if
$G \in {\cal OP}$. Thus it is interesting to decide whether the maximum possible value of the boxicity of
series-parallel graphs is $2$ or $3$. In this paper we construct a series-parallel graph with boxicity
$3$, thus resolving this question. Recently Chandran and Sivadasan \cite{sunns2} showed that for any $G$, 
$\mbox{box}(G) \le \mbox{treewidth}(G) + 2$. They conjecture that for any $k$, there exists a $k$-tree 
with boxicity $k+1$. (This would show that their upper bound is tight but for an additive factor of
$1$, since the treewidth of any $k$-tree equals $k$.) The series-parallel graph we construct in this paper
is a $2$-tree with boxicity $3$ and is thus a first step towards proving their conjecture. \\

\end{abstract}

\section{Introduction}

Let ${\cal F} = \brac{S_x \subseteq U: x \in V}$ be a family of subsets of a universe $U$, where
$V$ is an index set. The intersection graph $\Omega({\cal F})$ of ${\cal F}$ has $V$ as a vertex set,
and two distinct vertices $x$ and $y$ are adjacent if and only if $S_x \cap S_y \not= \emptyset$.
A $k$-dimensional box is a Cartesian product $R_1 \times R_2 \times \cdots \times R_k$ where
$R_i$ (for $1 \le i \le k$) is a closed interval of the form $\sbrac{a_i, b_i}$ on the real line.
For a graph $G$, its boxicity is the minimum dimension $k$, such that there exists a family 
${\cal F}$ of $k$-dimensional axis-parallel boxes with $\Omega({\cal F}) = G$.  We denote the boxicity of
a graph $G$ by box($G$). The notion of boxicity was introduced by Roberts 
\cite{Roberts} and has since
been studied by many authors such as Cozzens \cite{Coz}, Trotter \cite{Trotter} etc.

The complexity of finding the boxicity of a graph was shown to be NP-hard by 
Cozzens \cite{Coz}. This
was later improved by Yannakakis \cite{Yan1} and finally by Kratochvil 
\cite{Kratochvil} who showed that deciding whether the
boxicity of a graph is at most $2$ itself is NP-complete. 

Researchers have also tried to generalize or extend the concept of boxicity
in various ways. The poset boxicity \cite{TroWest}, the rectangle number 
\cite{ChangWest},
grid dimension \cite{Bellantoni}, circular dimension \cite{Feinberg,Shearer} and the boxicity of digraphs
 \cite{ChangWest1} are some examples.
\noindent
\vspace{0.32cm}\\
\textbf{Our Result} \\
Outer planar graphs (${\cal P}$), series-parallel graphs (${\cal SP})$ and planar graphs (${\cal OP}$) are 
three well-studied graph
classes. There is a natural hierarchy of proper inclusion relation among these graph
classes: ${\cal OP} \subset {\cal SP} \subset {\cal P}$. This
hierarchy looks quite natural as is evidenced by the well-known forbidden minor characterizations
of these graph classes: planar graphs are exactly the class of graphs with neither $K_5$ nor
$K_{3,3}$ as a minor; series-parallel graphs are exactly the graphs without a $K_4$ as a minor and 
outer planar graphs consists of exactly those graphs with neither $K_4$ nor $K_{3,2}$ as minor.
Two of the early results in the boxicity literature concern with the boxicity of planar
graphs and outer planar graphs. \\

\begin{theorem}
(Scheinerman \cite{Scheiner}).If $G$ is outer planar then box($G$) $\le 2$.
\end{theorem}
\begin{theorem}
\label{planar}
(Thomassen \cite{Thoma1}).If $G$ is a planar graph then box ($G$) $\le 3$.
\end{theorem}

Surprisingly, we haven't seen any attempts in the literature to decide whether the tightest
possible upper bound for the boxicity of series-parallel graphs is $2$ or $3$. Considering
the simple inductive definition of series-parallel graphs, (see Definition 1 below), one 
is tempted to believe that series-parallel graphs
 have
boxicity $2$, i.e. they can be represented as an intersection graph of axis-parallel rectangles.
Moreover, experimentation with small, easily constructible series-parallel graphs seems
to support this initial intuition. In this paper we construct a series-parallel graph
whose boxicity equals $3$. The series-parallel graph $G$ that we construct is fairly large 
(contains 157 vertices and 311 edges), $2$-connected and edge maximal.
In fact it is a $2$-tree. The reader may note that any series-parallel graph which
contains $G$ as an induced sub-graph also has boxicity $3$, and thus there exists an infinite family of 
series-parallel graphs with boxicity $3$.

The class of undirected graphs known as $k$-trees is defined recursively as follows: A $k$-tree
on $(k+1)$ vertices consists of a clique on $(k+1)$ vertices. Given any $k$-tree
$T_n$ on $n$ vertices ($n \ge k+1$) we construct a $k$-tree on $n+1$ vertices by adjoining a
new vertex $x_{n+1}$ to $T_n$, which is made adjacent to each vertex of some $k$-clique of
$T_n$ and non-adjacent to the remaining $n-k$ vertices.

It is well-known that the treewidth of a $k$-tree equals $k$ (see \cite{Bodland3} for a brief
survey on treewidth).
In fact the treewidth of a graph
$G$ can be defined as the smallest integer $k$, such that $G$ is a subgraph of some $k$-tree.
(A graph $G$ with treewidth $\le k$ is also known as a partial $k$-tree). Chandran and Sivadasan
 \cite{sunns2} have recently proved the following theorem.

\begin{theorem}
For any graph $G$, $\mbox{box}(G) \le \mbox{treewidth}(G) + 2$.
\end{theorem}

They construct a family of graphs such that $\mbox{box}(G)$ $\ge$ treewidth $(G)(1-o(1))$, establishing
the near-tightness of their result. On the other hand, they conjecture that their upper bound is tight but
for an additive factor of $1$: In particular they believe that for each $k \ge 1$, there exists a 
$k$-tree with boxicity $k+1$. The case $k=1$ is trivial since there are obviously $1$-trees 
(normal trees) whose boxicity equals $2$. It is well-known that the class of series-parallel graphs
is exactly the class of graphs with treewidth at most $2$. That is series-parallel graphs are
exactly the partial $2$-trees. Thus it is easy to see that every edge maximal series-parallel
graph is a $2$-tree. The graph we construct in this paper is an edge maximal series-parallel
graph and is thus a $2$-tree, whose boxicity equals $3$. Thus the construction given in this paper
settles the conjecture for $k=2$.

\begin{definition}
A connected series-parallel graph is a multigraph that can be constructed from a single vertex
by a sequence of applications of the following three operations
\begin{enumerate}
\item Series operation $($on an edge $(u,v))$: Add a new vertex $y$,
new edges $(u,y)$, $(y,v)$ and remove $(u,v)$.
\item Parallel operation $($on an edge $e=(u,v))$: Add another edge between $u \mbox{ and }v$.
\item Adding a pendant vertex $($to a vertex $u)$: Add a new vertex $y$ and a new edge
$(u,y)$.
\end{enumerate}
\end{definition}

\noindent
\textit{Remark:} Though series-parallel graphs are defined to be multigraphs,
their boxicity depends only on the underlying simple graphs.

\begin{definition}
$I = (V,E)$ is an interval graph if and only if there exists a function
$\Pi$ that maps each vertex $u \in V$ to a closed
interval of the form $[l(u), r(u)]$ on the real line such that $(u,v) \in E(I)
\iff \Pi(u) \cap \Pi(v) \not= \emptyset$. We will call $\Pi$,  an
interval representation of $I=(V,E)$.
\end{definition}

\begin{definition}
A $d$-box representation of $G=(V,E)$ is a function $\theta$ that  maps each vertex $v \in V(G)$ to a 
$d$-dimensional axis parallel box $R_1 \times R_2 \times \cdots R_d$, where $R_i$, for $1 \le i \le d$,
is a closed interval of the form $\sbrac{a_i, b_i}$ on the real line, such that $(u,v) \in E(G) \iff
\theta(u) \cap \theta(v) \not= \emptyset$. Let $\Pi_i$ be the function that maps $u \in V(G)$ to $R_i$.
Note that $\Pi_i(u)$ represents the projection of the box $\theta(u)$ on the $i$-th axis. The reader 
may also note that 
$\Pi_i$ is the interval representation of a graph $G' = (V,E')$ where $E' \supseteq E$. We write 
$\theta = (\Pi_1, \cdots , \Pi_d)$
\end{definition}

\begin{definition}
Boxicity of a graph $G$ is defined as the minimum $d$ such that a $d$-box representation exists for $G$, and is denoted by box$(G)$.
\end{definition}

\begin{definition}
A graph $G$ with box$(G) \le 2$ is called a rectangle graph. (This terminology is due to 
\cite{Trotter})
\end{definition}
A $2$-box representation of $G$ will also be called a rectangle representation of $G$. Since in this
paper we are dealing only with $2$-box representations, $\theta$ will always denote some rectangle
representation of $G$.

\begin{definition}
A \textit{split} operation on an edge $(a,b)$ of $G$ is the addition of a
new vertex $c$ to $V(G)$ and new edges $(a,c)$ and $(b,c)$ to $E(G)$.
We say that vertex $c$ is obtained by splitting $(a,b)$.
\end{definition}

Note that the \textit{split} operation on an edge $(a,b)$ is equivalent to a parallel operation on
$(a,b)$ followed by a series operation on the resulting (parallel) edge.

\section{Boxicity of Series-Parallel Graphs}
In this section, we construct a series-parallel graph with boxicity $ > 2$.
First we construct four graphs $L_1$ to $L_4$ which will occur as induced
subgraphs of the final graph.
Each graph $L_i$ has a bit larger size and has a bit more complex structure
than the previous graphs $L_j$, $j < i$.
As the graphs become more complex, we show that stricter constraints get imposed on their possible rectangle
representations. 
\noindent
\vspace{0.32cm}\\
\textbf{Construction of the graph $L_1$}: 
Start with an edge $(a,b)$, \textit{split} it to obtain a vertex $c$,
and add a pendant vertex $z$ to $c$. Thus $L_1$ has $4$ vertices and $4$ edges.

\begin{lemma}
\label{pendant}
Let $\theta$ be a rectangle representation of $L_1$. Then
$\theta(c) \not\subseteq \theta(a) \cup \theta(b)$.
\end{lemma}

\begin{proof}
Suppose $\theta(c) \subseteq \theta(a) \cup \theta(b)$.
Then $\theta(z) \cap \theta(c) \subseteq \theta(z) \cap (\theta(a) \cup \theta(b))$
$ = (\theta(z) \cap \theta(a))\cup(\theta(z) \cap \theta(b)) = \emptyset$ 
(since $(z,a),(z,b) \not\in E(L_1)$).
Thus $\theta(z) \cap \theta(c) = \emptyset$ which is a contradiction, since $(z,c) \in E(L_1)$.
\end{proof}

\begin{definition}
A family $\brac{T_i}_{i\in I}$ of subsets of a set $T$ has Helly Property if
for every $J \subseteq I$ the assumption that $T_i \cap T_j \not= \emptyset$
for every $i,j \in J$ implies ${\bigcap}_{j\in J}T_j \not= \emptyset$.
\end{definition}

It is easy to verify that a family of closed intervals on the real line satisfy Helly property. 
Now it is not difficult to infer that a family of $d$-dimensional axis parallel boxes also
satisfy Helly property. In particular, we have the following Lemma:

\begin{lemma}
\label{Helly}
Let $G$ be a triangle with vertices $a,b,c$. Let $\theta$ be a
rectangle representation of $G$. Then $\theta(a) \cap \theta(b) \cap \theta(c)
\not= \emptyset$.
\end{lemma} 
\noindent
\textbf{Construction of the graph $L_2$}: 
Start with an edge $(a,b)$, \textit{split} it to obtain a new vertex $c$, \textit{split} 
$(a,c)$ to obtain a new vertex $x$ and \textit{split} $(b,c)$ to obtain a new vertex $y$. 
The resulting graph $L_2$ has $5$ vertices and $7$ edges.
 
\begin{lemma}
\label{difference}
Let $\theta$ be a rectangle representation of $L_2$,
Then $\theta(c) \cap (\theta(a) - \theta(b)) \not= \emptyset$
and $\theta(c) \cap (\theta(b) - \theta(a)) \not= \emptyset$.
\end{lemma}

\begin{proof}
Suppose $\theta(c) \cap (\theta(a) - \theta(b)) = \emptyset$. Then we have
$\theta(c) \cap \theta(a) \subseteq \theta(b)$ and hence
$(\theta(c) \cap \theta(a))\cap \theta(x) \subseteq \theta(b) \cap \theta(x)$.
But $x,a,c$ induce a triangle in $L_2$ and therefore by Lemma \ref{Helly}, 
$\theta(x) \cap \theta(a) \cap \theta(c) \not= \emptyset$. On the other hand,
$\theta(x) \cap \theta(b) = \emptyset$ (since $(x,b) 
\notin E(L_2)$) which is a contradiction. 
Thus we can infer that $\theta(c) \cap(\theta(a) - \theta(b)) \not= \emptyset$.
Similarly we can show that  $\theta(c) \cap(\theta(b) - \theta(a)) \not= \emptyset$.
\end{proof}
\noindent
The following Lemma is intuitive. We prove it formally below.
\begin{lemma}
\label{projection}
Let $\theta = (\Pi_1, \Pi_2)$ be a rectangle representation of a graph $G$.
Then $\theta(c) \cap (\theta(a) - \theta(b)) \not= \emptyset$ if and only if 
at least one of the following two conditions hold 
\begin{enumerate}
\item $\Pi_1(c) \cap (\Pi_1(a) - \Pi_1(b)) \not= \emptyset$
\item $\Pi_2(c) \cap (\Pi_2(a) - \Pi_2(b)) \not= \emptyset$
\end{enumerate}
\end{lemma}

\begin{proof}
It is easy to verify that 
$\theta(c) \cap (\theta(a)- \theta(b)) \not= \emptyset$ if and only if
$\theta(c) \cap \theta(a) \not\subseteq \theta(b)$. This holds if and only if
$\Pi(a) \cap \Pi(c) \not\subseteq \Pi(b)$, for some $\Pi \in \brac{\Pi_1, \Pi_2}$,
since $\theta(u) = \Pi_1(u) \times \Pi_2(u)$ for any vertex $u$. But 
$\Pi(a) \cap \Pi(c) \not\subseteq \Pi(b)$ if and only if
$\Pi(c) \cap (\Pi(a) - \Pi(b)) \not= \emptyset$, and the Lemma follows.
\end{proof}

\begin{definition}
Let $\theta = (\Pi_1, \Pi_2)$ be a rectangle representation of a graph $G$, 
and let $a,b \in V(G)$ such
that $\theta(a) \cap \theta(b) \not= \emptyset$. Let $\Pi_1(a) = \sbrac{l_1(a), r_1(a)}$ ,
$\Pi_2(a) = \sbrac{l_2(a), r_2(a)}$,$\Pi_1(b) = \sbrac{l_1(b), r_1(b)}$ and 
$\Pi_2(b) = \sbrac{l_2(b), r_2(b)}$. Now if 
$l'_1 = max(l_1(a), l_1(b))$, $l'_2=max(l_2(a), l_2(b))$, 
$r'_1 = min(r_1(a), r_1(b))$ and $r'_2 = min(r_2(a), r_2(b))$. Then the corner points
of $\theta(a) \cap \theta(b)$ are defined to be the four points $(l'_1,l'_2),(l'_1, r'_2)
,(r'_1, l'_2),(r'_1, r'_2)$.
\end{definition}
Intuitively, the corner points are the four corners of $\theta(a) \cap \theta(b)$.
The four corner points need not be distinct.

\begin{lemma}
\label{corner}
Let $\theta = (\Pi_1, \Pi_2)$ be a rectangle representation of $L_2$.  
If $\Pi_i(c) \not\subseteq \Pi_i(a) \cap \Pi_i(b)$ for $ i = {1,2}$, then 
$\theta(c)$ contains a corner point of $\theta(a) \cap \theta(b)$.
\end{lemma}

\begin{proof}
First note that $(a,b) \in E(L_2)$ and therefore $\theta(a) \cap \theta(b) \not= \emptyset$ and
thus the corner points of $\theta(a) \cap \theta(b)$ are defined.
For $i = {1,2}$, let $l_i' = max(l_i(a), l_i(b))$ and $r_i' = min(r_i(a), r_i(b))$ represent the left and
right end points of $\Pi_i(a) \cap \Pi_i(b)$ respectively. In $L_2$, vertices $a,b,c$ 
induce a triangle and 
hence by Helly property (Lemma \ref{Helly}) we have $\theta(c) \cap \theta(a) \cap \theta(b) \not= \emptyset$.
It follows that for $i=1,2$,  $\Pi_i(c) \cap (\Pi_i(a) \cap \Pi_i(b)) \not= \emptyset$. Moreover, by assumption we have
for $i=1,2$, $\Pi_i(c) \not\subseteq \Pi_i(a) \cap \Pi_i(b)$.
Thus $\Pi_i(c)$ is an interval which contains at least one point from $\Pi_i(a) \cap \Pi_i(b)$ and
at least one point from the complement of $\Pi_i(a) \cap \Pi_i(b)$.
Therefore we can infer that either $l_i' \in \Pi_i(c)$ or $r_i' \in \Pi_i(c)$.
Thus we conclude that $\theta(c) = \Pi_1(c) \times \Pi_2(c)$ contains at least one
of the corner points $(l_1', l_2'), (l_1', r_2'), (r_1', l_2'), (r_1', r_2')$. 
\end{proof}
\noindent
\vspace{0.32cm}\\
\textbf{Construction of the graph $L_3$}: 
Start with a single edge $(a,b)$, \textit{split}
$(a,b)$ $5$ times to obtain the vertices $\brac{c_i: 1 \le i \le 5}$. For each
$c_i$, obtain $x_i$ by \textit{splitting} $(a,c_i)$ and $y_i$ by \textit{splitting}
$(b,c_i)$. Note that for $1 \le i \le 5$, $a,b,c_i,x_i,y_i$ induce a graph isomorphic to
$L_2$ in $L_3$.

\begin{definition}
Let $G=(V,E)$ be a graph with box$(G) \le 2$. Let
$\theta=(\Pi_1, \Pi_2)$ be a rectangle representation of $G$. We say that two vertices
$u,v \in V(G)$ are a \textit{crossing pair} with respect to $\theta$ if and only if one
of the following 2 conditions hold: 
\begin{enumerate}
\item $\Pi_1(u) \subseteq \Pi_1(v)$ and $\Pi_2(v) \subseteq \Pi_2(u)$ or
\item $\Pi_1(v) \subseteq \Pi_1(u)$ and $\Pi_2(u) \subseteq \Pi_2(v)$. 
\end{enumerate}
\end{definition}

\begin{lemma}
\label{cross}
Let $\theta=(\Pi_1, \Pi_2)$ be any rectangle representation of $L_3$. 
Then $a,b$ cannot be a \textit{crossing pair} with respect to $\theta$.
\end{lemma}

\begin{proof}
Suppose $a,b$ be a crossing pair. Then without loss of generality assume that 
\begin{eqnarray}
\Pi_1(a) \subseteq \Pi_1(b) \mbox{ and } \Pi_2(b) \subseteq \Pi_2(a)
\end{eqnarray}

\noindent
Now observe that for each $i$,  $1 \le i \le 5$, $a,b,c_i,x_i,y_i$ induce a subgraph isomorphic to $L_2$. 
Hence by Lemma \ref{difference}, we have 
\begin{eqnarray}
\theta(c_i) \cap (\theta(a) - \theta(b)) &\not=& \emptyset \mbox { and } \\
\theta(c_i) \cap (\theta(b) - \theta(a)) &\not=& \emptyset
\end{eqnarray}

\noindent
By Lemma \ref{projection} inequality (2) implies that at least one of the two
conditions (a) $\Pi_1(c_i) \cap (\Pi_1(a) - \Pi_1(b)) \not= \emptyset$ 
(b) $\Pi_2(c_i) \cap (\Pi_2(a) - \Pi_2(b)) \not= \emptyset$  holds.
But by condition (1), $\Pi_1(a) \subseteq \Pi_1(b)$, and hence 
$\Pi_1(c_i) \cap (\Pi_1(a) - \Pi_1(b)) = \emptyset$. Thus we infer that 
\begin{eqnarray}
\Pi_2(c_i) \cap (\Pi_2(a) - \Pi_2(b)) \not= \emptyset 
\end{eqnarray}
\noindent
In a similar way, from inequality (3) we can infer that
\begin{eqnarray}
\Pi_1(c_i) \cap (\Pi_1(b) - \Pi_1(a)) \not= \emptyset
\end{eqnarray}
\noindent
From inequalities (4) and (5) we get
\begin{eqnarray}
\Pi_j(c_i) \not\subseteq \Pi_j(a) \cap \Pi_j(b) \mbox{ for } j = 1,2
\end{eqnarray}
Therefore by Lemma \ref{corner}, for each $i$, $1 \le i \le 5$, $\theta(c_i)$ contains
a \textit{corner point} of $\theta(a) \cap \theta(b)$.
But since there are only at most $4$ corner points, by pegion hole principle there exist 
$i,j$ where $1 \le i,j \le 5$ and $i \not= j$ such that $\theta(c_i)$ and $\theta(c_j)$ contain
the same corner point, i.e. $\theta(c_i) \cap \theta(c_j) \not= \emptyset$. This
is a contradiction since $(c_i, c_j) \not\in E(L_3)$.
\end{proof} 
\noindent
\vspace{0.32cm}\\
\textbf{Construction of the graph $L_4$}: 
The graph $L_4$ is obtained from $L_3$ by splitting the edge $(x_i, c_i)$
to obtain $z_i$ for $1 \le i \le 5$.

\begin{lemma}
\label{main}
Let $\theta = (\Pi_1, \Pi_2)$ be a rectangle representation of $L_4$.
Then there exists $c \in \brac{c_i: 1 \le i \le 5}$ such that either $a,c$ 
or $b,c$ is a \textit{crossing pair}.
\end{lemma}

\begin{proof}
\noindent
We claim that there exists a $c \in \brac{c_i: 1 \le i \le 5}$ such that 
$\Pi_1(c) \subseteq \Pi_1(a) \cap \Pi_1(b)$ or  
$\Pi_2(c) \subseteq \Pi_2(a) \cap \Pi_2(b)$. Suppose not. Then 
for each $c_i$ and for $j=1,2$ , 
$\Pi_j(c_i) \not\subseteq \Pi_j(a) \cap \Pi_j(b)$. Thus by Lemma \ref{corner},
for each $i$, $1 \le i \le 5$, $\theta(c_i)$ contains a 
\textit{corner point} of $\theta(a) \cap \theta(b)$.
Since there are only at most $4$ corner points of $\theta(a) \cap \theta(b)$, by 
pegion hole principle,
there exist $i,j$ where $1\le i,j \le 5$ and $i \not=j$ such that 
$\theta(c_i) \cap \theta(c_j) \not= \emptyset$. This is a contradiction since
$(c_i, c_j) \not= E(L_4)$. 
Therefore without loss of generality we can assume that 
\begin{eqnarray}
\Pi_1(c_1) \subseteq \Pi_1(a)\cap \Pi_1(b)
\end{eqnarray}

\noindent
Now $\brac{a,b,c_1,x_1,y_1}$ induce a graph isomorphic to $L_2$ in $L_4$.
Therefore by Lemma \ref{difference}, $\theta(c_1) \cap (\theta(a) - \theta(b)) \not= \emptyset$ and
$\theta(c_1) \cap (\theta(b) - \theta(a)) \not= \emptyset$. By Lemma \ref{projection}, the
former inequality implies that at least one of the two conditions 
(a) $\Pi_1(c_1) \cap (\Pi_1(a) - \Pi_1(b)) \not= \emptyset$ 
(b) $\Pi_2(c_1) \cap (\Pi_2(a) - \Pi_2(b)) \not= \emptyset$ holds. But by condition (7),
$\Pi_1(c_1) \subseteq \Pi_1(a) \cap \Pi_1(b)$, and hence
$\Pi_1(c_1) \cap (\Pi_1(a) - \Pi_1(b)) = \emptyset$. Thus we infer that
\begin{eqnarray}
\Pi_2(c_1) \cap(\Pi_2(a) - \Pi_2(b)) \not= \emptyset
\end{eqnarray}

\noindent
Similarly $\theta(c_1) \cap (\theta(b) - \theta(a)) \not= \emptyset$ implies the following:
\begin{eqnarray}
\Pi_2(c_1) \cap (\Pi_2(b) - \Pi_2(a)) \not= \emptyset
\end{eqnarray}
\noindent
From (8) and (9) we claim that
\begin{eqnarray}
\Pi_2(a) \cap \Pi_2(b) \subseteq \Pi_2(c_1)
\end{eqnarray}
\noindent
To verify the above, let $l'$ and $r'$ be the left and right endpoints respectively of 
$\Pi_2(a) \cap \Pi_2(b)$. ($\Pi_2(a) \cap \Pi_2(b) \not= \emptyset$ since
$\theta(a) \cap \theta(b) \not= \emptyset$). Let $x \in (\Pi_2(a) - \Pi_2(b)) \cap \Pi_2(c_1)$
and $y \in (\Pi_2(b) - \Pi_2(a)) \cap \Pi_2(c_1)$. (Inequalities (8) and (9) ensure that
we can find such an $x$ and $y$).  Since $\Pi_2(a)$ and $\Pi_2(b)$
are intervals it is easy to verify that either $x < l' < r' < y$ or 
$y < l' < r' < x$. Without loss of generality let $x < l' < r' < y$. Then, clearly we have,
$\Pi_2(c_1) \supseteq \sbrac{x,y} \supset \sbrac{l',r'} = \Pi_2(a) \cap \Pi_2(b)$, as required.

\vspace{0.32cm}
\noindent
Observe that the graph induced by $\brac{a,b,c_1,z_1}$ in $L_4$ is isomorphic to $L_1$.
Hence by Lemma \ref{pendant}, $\theta(c_1) \not\subseteq \theta(a) \cup \theta(b)$. 
Since by (7), $\Pi_1(c_1) \subseteq \Pi_1(a) \cap \Pi_1(b)$ we must have 
\begin{eqnarray}
\Pi_2(c_1) \not\subseteq \Pi_2(a) \cup \Pi_2(b)
\end{eqnarray}

\noindent
Let $l'',r''$ be the left and right end points respectively
of $\Pi_2(a) \cup \Pi_2(b)$. It is easy to see that the set $\brac{l',r',l'',r''}$ is the same as the set 
$\brac{l_2(a), l_2(b), r_2(a), r_2(b)}$. Since $\Pi_2(a), \Pi_2(b)$
and $\Pi_2(c_1)$ are intervals, (10) and (11) imply that
at least $3$ of these points are contained in $\Pi_2(c_1)$. Thus either 
$\sbrac{l_2(a), r_2(a)} \subseteq \Pi_2(c_1)$ or $\sbrac{l_2(b), r_2(b)} \subseteq \Pi_2(c_1)$.
In other words:
\begin{eqnarray}
\Pi_2(a) \subseteq \Pi_2(c_1)\mbox{ or }\Pi_2(b) \subseteq \Pi_2(c_1)
\end{eqnarray}
\noindent
By (7) and (12), we conclude that either $a,c_1$ is a \textit{crossing pair} or 
$b,c_1$ is a \textit{crossing pair}.

\end{proof}

\subsection{A Series-Parallel graph whose boxicity $> 2$.}
Using the ideas presented above we present a series-parallel graph whose
boxicity $>$ 2. The construction is as follows.

\begin{enumerate}
\item Let the initial graph be the single edge $(a,b)$.
\item For $i = 1 \mbox { to } 5$ do:\\
Apply the \textit{split} operation on $(a,b)$ and let $c_i$ be the resulting vertex. 
\item For each $c_i$ where $1 \le i \le 5$ do 
\begin{enumerate}
\item Apply the \textit{split} operation on $(a,c_i)$ five times: Let $d_{ij}$ where $1 \le j \le 5$
be the resulting vertices.
\item Apply the \textit{split} operation on $(b,c_i)$ five times: Let $e_{ij}$ where $1 \le j \le 5$
be the resulting vertices.
\end{enumerate}

\item For all $i,j$ where $1 \le i,j \le 5$ \\
Apply the \textit{split} operation on $(a, d_{ij})$, $(c_i, d_{ij})$, $(b, e_{ij})$ and
$(c_i, e_{ij})$. Let the resulting vertices be $p_{ij}$, $q_{ij}$, $r_{ij}$ and $s_{ij}$ respectively.
\end{enumerate}

Note that the graph $G$ constructed above is a series-parallel graph, since we are using the
\textit{split} operations only. Also note that $G$ has $n=157$ vertices and
$2n-3=311$ edges. Since any series-parallel graph on $n$ vertices
with $2n-3$ edges is edge maximal (see chapter 8, Diestel \cite{Diest}), it follows that 
$G$ is an edge maximal series-parallel graph. Thus it is also a $2$-tree (This fact 
is in fact evident from the construction since we are using \textit{split} operations only) 
and hence $2$-connected.

\begin{theorem}
The graph $G$ defined above has boxicity $ = 3$.
\end{theorem}
\begin{proof}
First we show that box$(G) > 2$. Suppose not. Then there exists a rectangle
representation for $G$.
It is easy to verify that $\brac{a,b} \cup \brac{c_i,d_{i1}, e_{i1}, q_{i1}: 1 \le i \le 5}$ 
induce a graph isomorphic to $L_4$.
Therefore by Lemma \ref{main}, there exists a $c \in \brac{c_i:1 \le i \le 5}$ such that either $a,c$
is
a \textit{crossing pair} or $b,c$ is a \textit{crossing pair}. Without loss
of generality let $a,c_1$ be a \textit{crossing pair}.
But $\brac{a,c_1} \cup \brac{d_{1j}, p_{1j}, q_{1j}: 1 \le j \le 5}$, induce a graph isomorphic
to $L_3$.  Thus by Lemma \ref{cross}, $a,c_1$ cannot be a \textit{crossing pair}, which is a contradiction.  
Thus we infer that box$(G) > 2$. 
Since any series-parallel graph is planar we have (by Theorem \ref{planar}) box$(G) \le 3$ 
and the theorem follows.
\end{proof}

\section{Conclusions and Open Problems}
In this paper we have shown that there exists an infinite family of series-parallel graphs
with boxicity equal to $3$. Thus the following problem arises naturally.

\noindent
1. Characterize the class of series-parallel graphs with boxicity $\le 2$. 

It is implicit in a Theorem of Thomassen \cite{Thoma1} that any series-parallel graph without the 
join of $K_2$ and $\bar{K_3}$ as an induced subgraph has a strict box representation. 
Another interesting open problem is:

\noindent
2. Prove that for each $k \ge 1$, there exist a $k$-tree with boxicity $=k+1$. 

The case $k=1$ is trivial, and the case $k=2$ is settled in this paper. Also it is not
difficult to show that there exist $k-$trees with boxicity at least $\floor{k/2}$.


\begin{thebibliography}{10}

\bibitem{Bellantoni}
S.~Bellantoni, I.~Ben-Arroyo Hartman, T.~Przytycka, and S.~Whitesides.
\newblock Grid intersection graphs and boxicity.
\newblock {\em Discrete mathematics}, 114(1-3):41--49, April 1993.

\bibitem{Bodland3}
H.~L. Bodlaender.
\newblock A tourist guide through treewidth.
\newblock {\em Acta Cybernetica}, 11:1--21, 1993.

\bibitem{sunns2}
L.Sunil Chandran and Naveen Sivadasan.
\newblock Treewidth and boxicity.
\newblock Preprint, 2005.

\bibitem{ChangWest1}
Y.~W. Chang and Douglas~B. West.
\newblock Interval number and boxicity of digraphs.
\newblock In {\em Proceedings of the 8th International Graph Theory Conf.
  (Kalamazoo 1996) (Wiley, 1998)}.

\bibitem{ChangWest}
Y.~W. Chang and Douglas~B. West.
\newblock Rectangle number for hyper cubes and complete multipartite graphs.
\newblock In {\em 29th {SE} conf. Comb., Graph Th. and Comp., Congr. Numer.
  132(1998)}, 19--28.

\bibitem{Coz}
M.~B. Cozzens.
\newblock Higher and multidimensional analogues of interval graphs.
\newblock Ph. D thesis, Rutgers University, New Brunswick, NJ, 1981.

\bibitem{Diest}
Reinhard Diestel.
\newblock {\em Graph Theory}, volume 173.
\newblock Springer Verlag, New York, 2 edition, 2000.

\bibitem{Feinberg}
Robert~B. Feinberg.
\newblock The circular dimension of a graph.
\newblock {\em Discrete mathematics}, 25(1):27--31, 1979.

\bibitem{Kratochvil}
J.~Kratochvil.
\newblock A special planar satisfiability problem and a consequence of its
  {NP}--completeness.
\newblock {\em Discrete Applied Mathematics}, 52:233--252, 1994.

\bibitem{Roberts}
F.~S. Roberts.
\newblock {\em Recent Progresses in Combinatorics}, chapter On the boxicity and
  Cubicity of a graph, pages 301--310.
\newblock Academic Press, New York, 1969.

\bibitem{Scheiner}
E.~R. Scheinerman.
\newblock Intersectin classes and multiple intersection parameters.
\newblock Ph. D thesis, Princeton University, 1984.

\bibitem{Shearer}
J.~B. Shearer.
\newblock A note on circular dimension.
\newblock {\em Discrete mathematics}, 29(1):103--103, 1980.

\bibitem{Thoma1}
C.~Thomassen.
\newblock Interval representations of planar graphs.
\newblock {\em Journal of combinatorial theory, Ser {B}}, 40:9--20, 1986.

\bibitem{TroWest}
Jr. W.~T.~Trotter and Douglas~B. West.
\newblock Poset boxicity of graphs.
\newblock {\em Discrete Mathematics}, 64(1):105--107, March 1987.

\bibitem{Trotter}
Jr~William~T.Trotter.
\newblock A characterization of robert's inequality for boxicity.
\newblock {\em Discrete Mathematics}, 28:303--313, 1979.

\bibitem{Yan1}
Mihalis Yannakakis.
\newblock The complexity of the partial order dimension problem.
\newblock {\em {SIAM} Journal on Algebraic Discrete Methods}, 3:351--358, 1982.

\end{thebibliography}
\end{document}